\documentclass{amsart}

\usepackage{amsmath}
\usepackage{amscd}
\usepackage{amssymb} 

\newcommand{\cal}{\mathcal}
\newcommand{\bk}{{\bf k}}

\newcommand{\bC}{{\Bbb C}}

\newcommand{\bQ}{{\Bbb Q}}
\newcommand{\bR}{{\Bbb R}}

\newcommand{\cA}{{\cal A}}
\newcommand{\cB}{{\cal B}}

\newcommand{\cH}{{\cal H}}
\newcommand{\cL}{{\cal L}}
\newcommand{\cM}{{\cal M}}
\newcommand{\cO}{{\cal O}}

\newcommand{\fd}{{\frak d}}

\newcommand{\fsl}{{\frak s}{\frak l}}

\DeclareMathOperator{\Img}{Im}

\DeclareMathOperator{\Ker}{Ker}

\newtheorem{theorem}{Theorem}[section]
\newtheorem{theorem/definition}{Theorem/Definition}[section]

\newtheorem{lemma}{Lemma}[section]

\newtheorem{corollary}{Corollary}[section]

\newenvironment{remark}{\medskip 
\noindent {\bf Remark.}}{\mbox{}}

 \newenvironment{example}{\medskip 
\noindent {\bf Example.}}{\mbox{}} 
\newenvironment{definition}{\medskip 
\noindent {\bf Definition.}}{\mbox{}} 

\begin{document}

\title
{Rational homotopy types of mirror manifolds}
\author{Jian Zhou}
\address{Department of Mathematics\\
Texas A\&M University\\
College Station, TX 77843}
\email{zhou@math.tamu.edu}
\begin{abstract}
We explain how to relate the problem of finding a mirror manifold
for a Calabi-Yau manifold to the problem of characterizing 
the rational homotopy types of closed K\"{a}hler manifolds.
\end{abstract}
\maketitle
\date{}
%\footnotetext[1]{1991 {\em Mathematics Subject Classification}: 53C15, 58A12,
%81R05}
%\footnotetext[1]{Author's research was supported in part by NSF}

In this paper, 
we show that under a rationality condition on the (classical)
Yukawa coupling constants for a Calabi-Yau manifold $M$, 
a rational homotopy type can be determined.
If $M$ has a mirror manifold $\widetilde{M}$,
then $\widetilde{M}$ has this rational homotopy type.
Thus the problem of 
finding a mirror manifold of a Calabi-Yau manifold 
is related to the problem of
characterizing the rational homotopy types of K\"{a}hler manifolds. 
 
We refer to Yau \cite{Yau} and Greene-Yau \cite{Gre-Yau} 
for the history and background on the mirror symmetry
suggested by string theory.
Mathematically,
there are two natural algebras defined by any 
Calabi-Yau $n$-fold,
\begin{eqnarray*}
A(M) = \oplus_{p, q \geq 0} H^q(M, \Omega^p), &
B(M) =  \oplus_{p, q \geq 0} H^q(M, \Omega^{-p}),
\end{eqnarray*}
where  $\Omega^{-p}$ is the sheaf of holomorphic sections
to $\Lambda^p TM$.
A weak version of mirror symmetry states 
that certain Calabi-Yau manifolds come in ``mirror pairs" $(M, \widetilde{M})$,
such that $A(M) \cong B(\widetilde{M})$ and $B(M) \cong A(\widetilde{M})$
as algebras.
String theory suggests the study of some geometrically defined
deformations of the algebra structures on $A(M)$ and $B(M)$
with prepotential functions
which satisfy the WDVV equations. 
(The notion of a Frobenius manifold introduced by Dubrovin \cite{Dub}
gives a coordinate free formulation of such deformations.)
On $A(M)$,
the quantum cohomology theory 
provides a construction of a structure of formal Frobenius manifold.
See e.g. Ruan-Tian \cite{Rua-Tia}.
For $B(M)$,
earlier study of mirror symmetry 
studies the variation of Hodge structures under the deformation 
of complex structure on a Calabi-Yau manifold,
using the Bogomolov-Tian-Todorov theorem \cite{Tia, Tod}.
See e.g. Morrison \cite{Mor}.
A famous prediction based on the conjectured mirror symmetry is 
the formula of Candelas {\em et al} \cite{Can-Oss-Gre-Par}
which has been rigorously proved by
Lian, Liu and Yau \cite{Lia-Liu-Yau}. 
Based on a key formula in Tian \cite{Tia},
Bershadsky-Cecotti-Ooguri-Vafa \cite{Ber-Cec-Oog-Vaf}
wrote a Lagrangian whose Euler-Lagrange equation is 
the deformation equation of the complex structure on a Calabi-Yau manifold.
They also noticed that this Lagrangian when restricted to 
the critical submanifold,
gives the prepotential function for a deformation of $B(M)$.
Mathematical formulation of this result together with
a version of mirror symmetry,
both in terms of Frobenius manifolds,
can be found in Barannikov-Kontsevich \cite{Bar-Kon}
where they emphasized on the role of the DGBV algebra structure revealed by
Tian's formula and remarked that such a construction can also be carried
out for DGBV algebras with suitable conditions.
For details,
see Manin \cite{Man2}.
More recently, 
 Cao and the author
\cite{Cao-Zho1, Cao-Zho2} observed that for any closed K\"{a}hler manifold $M$,
$(\Omega^{*, *}(M), \wedge, \bar{\partial})$ 
can be given a structure of a DGBV algebra
which satisfies the conditions to carry out the construction of a
formal Frobenius manifold structure on the Dolbeault cohomology
$A(M) = H^{*, *}(M)$.
See also Cao-Zhou \cite{Cao-Zho3}. 
For Calabi-Yau manifolds,
we conjecture that this Frobenius manifold structure 
can be identified with that provided by quantum cohomology. 
Also,
we formulate the following version of mirror symmetry : 
the Frobenius manifold constructed 
by Barannikov and Kontsevich \cite{Bar-Kon} for a Calabi-Yau manifold
can be identified 
with that constructed in \cite{Cao-Zho1, Cao-Zho2} for its mirror manifold.

Our version of mirror symmetry conjecture
can be studied by techniques in the rational homotopy theory.
Cao and the author \cite{Cao-Zho4} introduced 
a natural notion of quasi-isomorphisms of DGBV algebras
and proved that quasi-isomorphic DGBV algebras give rise to
isomorphic formal Frobenius manifold structures.
We also used Chen's approach to rational homotopy theory 
to formulate generalized version of mirror symmetry between complex manifolds 
and symplectic manifolds \cite{Zho}.
Such progress strongly motivates this work.
We show that under some rationality condition,
the cohomology algebra $B(M)$ determines a rational homotopy type.
If there is a K\"{a}hler manifold with this rational homotopy type,
then it is a mirror manifold of $M$;
on the other hand,
if $M$ has a mirror manifold $\widetilde{M}$,
then $\widetilde{M}$ has this rational homotopy type.
A new perspective to the topology change phenomenon
discussed in Aspinwall-Greene-Morrison \cite{Asp-Gre-Mor} and Witten \cite{Wit}
is that though the differential structures of the Calabi-Yau manifolds involved
may change,
the rational homotopy types probably remains the same.
Our result does not tell how to construct the mirror manifold,
since it is not known what kind of rational homotopy types are realized by
K\"{a}hler manifolds. 
Indeed an unsettled question of Sullivan asks for the classification of
rational homotopy types realized by  closed symplectic manifolds.

The rest of the paper is arranged as follows.
We review in \S 1 
some basics of the rational homotopy theory,
in particular, the notion of a formal differential graded algebra (DGA)
is recalled.
A method of proving formality  due to 
Deligne-Griffiths-Morgan-Sullivan \cite{Del-Gri-Mor-Sul} is 
formulated in a general situation in \S 2.
\S 3 contains the applications to K\"{a}hler and Calabi-Yau manifolds.
We give some postulations on the topology change from
consideration of complex and K\"{a}hler moduli space of a Calabi-Yau
manifold in \S 4.

{\bf Acknowledgements}. 
{\em The author thanks Texas A$\&$M University
and the Department of Mathematics there
for hospitality and financial support.
This research is partly supported by a NSF group infrastructure grant
through the Geometry-Analysis-Topology group at TAMU.
Special thanks are due to Huai-Dong Cao for a wonderful 
experience of collaboration which leads to this work.
The author also thanks him for suggestions which improves the 
presentation of this paper.}

\section{A brief review of rational homotopy theory}

In this section, we review some 
definitions and well-known results in rational homotopy theory.
The original sources  are Quillen \cite{Qui},
Sullivan \cite{Sul}, 
Deligne-Griffiths-Morgan-Sullivan \cite{Del-Gri-Mor-Sul}, 
Griffiths-Morgan \cite{Gri-Mor} and Wu \cite{Wu}.

\begin{definition}
A {\em differential graded algebra} (DGA) is a graded vector
space over $\bk = \bQ$, $\bR$ or $\bC$, 
$$ \cA = \oplus_{p \geq 0} \cA^p,$$
together with a linear operator
$d: \cA \rightarrow \cA$ of degree $1$ with $d^2 = 0$, 
a product $\wedge: \cA^p \otimes \cA^q \rightarrow \cA^{p+q}$,
such that
\begin{align*}
&\alpha \wedge \beta = (-1)^{|\alpha||\beta|} \beta \wedge \alpha, 
& \text{(graded commutativity)} \\
&d (\alpha \wedge \beta) = d\alpha \wedge \beta 
+ (-1)^{|\alpha|} \alpha \wedge d \beta. &
\text{(derivation property)}
\end{align*}
The cohomology of $(\cA, \wedge, d)$ is defined as usual by
$$H^*(\cA) : = \Ker d / \Img d.$$
The multiplication $\wedge$ on $\cA$ induces a graded commutative 
multiplication on $H^*(\cA)$.
If $H^0(\cA) = \bk$, then $\cA$ is called {\em connected}.
If we also have $H^1(\cA) = 0$, 
then $\cA$ is called {\em simply connected}.
We omit the routine definitions of DGA-homomorphisms,
DGA-isomorphisms and tensor products of DGA's. 
\end{definition}

\begin{example} (i) The real or complex de Rham algebra of a smooth 
manifold is a DGA with the exterior product and the exterior differential.

(ii) Any graded commutative algebra can be regarded as a DGA with
$d = 0$, e.g. the cohomology $H^*(\cA)$ with the induced multiplication
for any DGA $(\cA, \wedge, d)$.
\end{example}

\begin{definition}
A DGA $\cM$ is called {\em minimal} if

(i) $\cM$ is free as a graded commutative algebra.

(ii) $\cM^1 = 0$, and

(iii) $d \cM \subset \cM^+ \wedge \cM^+$,
where  $\cM^+ = \oplus_{p > 0} \cM^p$.
\end{definition}

Condition (i) means that $\cM$ is a tensor product of the polynomial 
algebra on some generators of even degrees with an exterior algebra
on some generators of odd degrees.
These generators span a linear subspace $I(\cM)$, 
called the {\em space of indecomposable elements} of $\cM$.
For each integer $k$, 
denote by $(\cM^+)^k$ the $k$-th power of $\cM^+$, 
i.e. the space spanned by products of $k$ elements in $\cM^+$. 
Then it is clear that we have isomorphisms
$$(\cM^+)^{k}/(\cM^+)^{k+1} \cong I(\cM) \wedge \cdots \wedge I(\cM).  
\,\,\,\,\,\, \text{($k$ times)}
$$
Condition (iii) implies that $d (\cM^+)^{k} \subset (\cM^+)^{k+1}$. 
Hence one gets an induced map 
$\bar{d}: I(M) \rightarrow I(M) \wedge I(M)$.
The couple $(I(\cM), \bar{d})$ is in fact an {\em shifted Lie coalgebra}, 
i.e. 
$[\cdot, \cdot] = \bar{d}^*: I(\cM)^* \otimes I(\cM)^* \rightarrow I(\cM)^*$
satisfies the graded Jacobi identity if the grading
of $I(\cM)$ is shifted by $1$.

\begin{definition}
The {\em minimal model} of a DGA $\cA$ is a DGA-homomorphism 
$\rho: \cM \rightarrow \cA$, such that $\cM$ is minimal,
and $\rho$  induces isomorphism on cohomology.
\end{definition}

A basic result in Sullivan's minimal model theory is the following

\begin{theorem}
Every simply connected DGA has a minimal model,
unique up to DGA-isomorphisms.
\end{theorem}

\begin{definition}
A DGA $\cA$ is called {\em formal} if it has the same minimal model 
as its cohomology algebra.
\end{definition}

\begin{definition}
A {\em quasi-isomorphism} between two DGA's $\cA$ and $\cB$ is a 
series of DGA's $\cA_0, \cdots, \cA_n$, 
and DGA-homomorphisms either $f_i: \cA_i \rightarrow \cA_{i+1}$, 
or $f_i: \cA_{i+1} \rightarrow \cA_i$ for $0 \leq i \leq n-1$,
such that $\cA_0 = \cA$, $\cA_n = \cB$ and each $f_i$ induces 
isomorphism on cohomology.
\end{definition}

Clearly, a formal DGA is quasi-isomorphic to its cohomology algebra.

\begin{theorem}
Quasi-isomorphic DGA's have the same minimal model.
\end{theorem}

\begin{corollary}
If a DGA is quasi-isomorphic to its cohomology algebra,
then it is formal.
\end{corollary}

\begin{example}
A well-known result in Deligne {\em et al} \cite{Del-Gri-Mor-Sul}
shows the de Rham algebra of a closed K\"{a}hler manifold is formal.
\end{example}

The homotopy groups of a topological space are usually very hard to compute.
For simplicity, we assume the space is simply connected. 
The work of Serre \cite{Ser} shows that it is easier to compute
the rank of homotopy groups by tensoring with $\bQ$.
Let
\begin{eqnarray*}
&& \pi_*(M) = \oplus_{p > 1} \pi_p(M), \\
&& \bQ \pi_*(M) = \oplus_{p > 1} \pi_p(M) \otimes \bQ.
\end{eqnarray*}
There is a well-known Whitehead product 
$$ [\cdot, \cdot]: \pi_p(M) \otimes \pi_q(M) \rightarrow \pi_{p+q-1}(M),$$
which extends naturally to $\bQ \pi_*(M)$.
The anti-symmetry and graded Jacobi identity for the Whitehead product
imply that $(\bQ \pi_*(M), [\cdot, \cdot])$ 
is an odd graded Lie algebra over $\bQ$.
We call it the {\em rational homotopy algebra} of $M$.
Similarly, one can define the {\em real homotopy algebra}
of $M$.
The following theorem is due to Serre \cite{Ser}:

\begin{theorem}
The following assertions are equivalent 
in the category of simply connected pointed topological space 
and continuous basepoint preserving maps:

(i) $\pi_*(f) \otimes \bQ: \bQ \pi_*(X) \rightarrow \bQ \pi_*(Y)$
is an isomorphism.

(ii) $H_*(f, \bQ): H_*(X, \bQ) \rightarrow H_*(Y, \bQ)$ is an
isomorphism.
\end{theorem}

A map satisfying these conditions are called a 
{\em rational homotopy equivalence},
and the spaces $X$ and $Y$ are then said 
to {\em have the same rational homotopy type}.
Equivalently, one can define the localization $X \rightarrow X(0)$
of a space $X$ at $0$. Then rational homotopy type of $X$ is 
the homotopy type of $X(0)$.
We say that an odd graded Lie algebra $\cL = \oplus_{p > 0} \cL_p$
is {\em reduced} if $\cL_1 = 0$.
The following important result is due to Quillen \cite{Qui} 
(p. 210, Corollary):

\begin{theorem}
Any reduced rational odd graded Lie algebra 
is the rational homotopy algebra 
of some simply connected pointed space $X$.
\end{theorem}

The main result of Sullivan's theory is that for a simply connected 
smooth manifold $M$,
the $\bQ$-polynomial de Rham algebra determines 
the rational homotopy type of $M$.
In particular, 
if  the $\bQ$-polynomial de Rham algebra is formal,
then the cohomology ring $H^*(M, \bQ)$ determines the rational
homotopy type of $M$.
An important result in Sullivan \cite{Sul} is that the formality
of the $\bQ$-polynomial de Rham algebra over $\bQ$ is equivalent to 
the formality of the de Rham algebra over $\bR$ or $\bC$,
at least for simply connected manifolds.

\section{A method of proving formality}

A well-known result in rational homotopy theory is 
that the de Rham algebra of 
a closed K\"{a}hler manifold is formal.
Deligne {\em et al}. \cite{Del-Gri-Mor-Sul} gave two proofs of this fact.
We recall in this section the method of their first proof
in a slightly more general situation.

Suppose that a differential graded algebra $(\cA, \wedge, \fd)$ is given
a (Euclidean or Hermitian) metric $\langle \cdot, \cdot \rangle$,
such that $\fd$ has a formal adjoint $\fd^*$,
i.e.,
$$\langle \fd a, b \rangle = \langle a, \fd^* b\rangle.$$
Set $\square_{\fd} = \fd\fd^* + \fd^* \fd$,
and $\cH = \Ker \square_{\fd} = \{a \in \cA: \square_{\fd} a = 0 \}$.
Then $\cH = \{a: \fd a = 0, \fd^* a = 0 \} = \Ker \fd \cap \Ker \fd^*$.
Assume that $\cA$ admits a ``Hodge decomposition'':
\begin{eqnarray} \label{eqn:HodgeDec1}
\cA = \cH \oplus \Img \fd \oplus \Img \fd^*.
\end{eqnarray}
It is standard to see that $\Ker \fd = \cH \oplus \Img \fd$,
and hence $H(\cA, \fd) \cong \cH$ as vector spaces.
For any $\alpha \in A$,
denote by $\alpha^H$ its projection onto $\cH$. 
For $\alpha, \beta \in \cH$, define
$$\alpha \circ \beta = (\alpha \wedge \beta)^H.$$
Note that there is an induced wedge product $\wedge$ on $H(A, \fd)$.
We have the following

\begin{lemma} \label{lm:isomorphism}
The isomorphism $\phi: \cH \to H(A, \fd)$ given by $\alpha \mapsto [\alpha]$
(the cohomology class of $\alpha$) maps $\circ$ to $\wedge$,
hence $(\cH, \circ)$ is an algebra which is isomorphic
to $(H(\cA, \fd), \wedge)$.
\end{lemma}

\begin{proof}
For $\alpha, \beta \in \cH$,
we have
\begin{eqnarray*}
\phi( \alpha \circ \beta) = [(\alpha \wedge \beta)^H] 
= [\alpha \wedge \beta] = [\alpha] \wedge [\beta] 
= \phi(\alpha) \wedge \phi(\beta).
\end{eqnarray*}
\end{proof}

Now assume that there is another differential $\fd^c$ on $(\cA, \wedge)$
with formal adjoint $(\fd^c)^*$,
such that the following identity holds:
\begin{eqnarray}
&& [\fd, \fd^c] = [\fd, (\fd^c)^*] = 0, \\
&& \square_{\fd} = \square_{\fd^c}.
\end{eqnarray}
By taking the formal adjoints, 
we also have
\begin{eqnarray}
&& [\fd^*, (\fd^c)^*] = [\fd^*, \fd^c] = 0.
\end{eqnarray}
Furthermore,
we assume that there is a ``Hodge decomposition" also for $d^c$:
\begin{eqnarray} \label{eqn:HodgeDec2}
\cA = \cH \oplus \Img \fd^c \oplus \Img (d^c)^*.
\end{eqnarray}
Combining with (\ref{eqn:HodgeDec1}),
we get a fivefold decomposition
\begin{eqnarray} \label{eqn:fivefold}
\cA = \cH \oplus \Img \fd\fd^c \oplus \Img \fd (\fd^c)^*
\oplus \Img \fd^*\fd^c \oplus \Img \fd^*(\fd^c)^*.
\end{eqnarray}

\begin{theorem} \label{thm:formal}
Under the above assumptions,
$(\cA, \wedge, \fd)$ is formal. 
\end{theorem}
 
\begin{proof} Since $\fd$ commutes with $\fd^c$,
$(\Ker \fd^c, \wedge,  \fd)$ is a DGA.
Consider the sequence
\begin{eqnarray} \label{eqn:quasi}
(\cA, \wedge, \fd)
\stackrel{i}{\hookleftarrow} (\Ker \fd^c, \wedge, \fd) 
\stackrel{\pi}{\rightarrow} (H^*(\cA, \fd^c), \wedge, 0),
\end{eqnarray}
where $i$ is the inclusion,
and $\pi$ is the quotient map.
Both $i$ and $\pi$ are homomorphisms of DGA's. 
It suffices to show that $\fd$ induces 
the zero homomorphism on $H^*(\cA, \fd^c)$.
This is easy to see 
since every element of $H^*(\cA, \fd^c)$ is represented by an element in $\cH$.
Furthermore, 
they are quasi-isomorphisms since all the cohomology groups
can be identified with $\cH$.
It is trivial for the first and third DGA's.
For the second DGA,
it follows from the decomposition
$$\Ker \fd^c = \cH \oplus \Img \fd\fd^c \oplus \Img \fd^*\fd^c.$$
By Lemma \ref{lm:isomorphism},
$(H^*(\cA, \fd^c), \wedge) \cong (\cH, \circ) \cong (H^*(\cA, \fd), \wedge)$.
This completes the the proof.
\end{proof}
 
\begin{remark}
Under the same assumptions,
$(\cA, \wedge, \fd^c)$ is also formal.
Furthermore,
by Lemma \ref{lm:isomorphism},
it has the same minimal model as $(\cA, \wedge, \fd)$.
\end{remark}

\section{Applications to K\"{a}hler and Calabi-Yau geometries}

\subsection{DGA $\cA$ of a K\"{a}hler manifold}

Let $M$ be a complex manifold,
and
$$\cA = \Omega^{*, *}(M) = \oplus_{p, q} \Omega^{p, q}(M) 
= \oplus_{p, q} \Gamma(M, \Lambda^pT^*M 
\otimes \Lambda^q\overline{T}^*M),$$
with ordinary wedge product $\wedge$.
Operators $d$ and $\bar{\partial}$ are differentials on $(\cA, \wedge)$.
As proved by Deligne {\em et al.} \cite{Del-Gri-Mor-Sul},
when $M$ is a closed K\"{a}hler manifold,
the DGA's thus obtained are formal.
One can make use the following 
well-known commutation relations in Hodge theory of K\"{a}hler manifolds
(see e.g. Griffiths-Harris \cite{Gri-Har}):
\begin{align}
& [\partial, \partial] = 0, & [\bar{\partial}, \bar{\partial}] &= 0, 
	& [\partial, \bar{\partial}] &= 0, 
	\label{eqn:comm1} \\
& [\partial^*, \partial^*] = 0, & [\bar{\partial}^*, \bar{\partial}^*] &= 0, 
	& [\partial^*, \bar{\partial}^*] &= 0, \\
& [\partial, \bar{\partial}^*] = 0, & [\partial^*, \bar{\partial}] &= 0,
	& \\
& [\partial, \partial^*] = \frac{1}{2}\square, 
& [\bar{\partial}, \bar{\partial}^*] & = \frac{1}{2}\square, 
	& \label{eqn:comm4}
\end{align}
where $\square$ is the Laplacian.
There are similar formulas for $d$ and $d^c$.
They imply that $\fd = \partial$ and $\fd^c = \bar{\partial}$ 
or $\fd = d$ and $\fd^c = d^c$ 
satisfy the conditions in Theorem \ref{thm:formal}.
So we have

\begin{theorem} \label{thm:formal2}
On a closed K\"{a}hler manifold, 
the DGA's $(\cA = \Omega^{*,*}(M), \wedge, \bar{\partial})$, 
$(\cA, \wedge, \partial)$,  
$(\cA, \wedge, d)$  and $(\cA, \wedge, d^c)$ are all formal.
Furthermore, they all have the same minimal model
as the Dolbeault cohomology algebra $(H^{*,*}(M), \wedge, 0)$.
\end{theorem}

The Dolbeault cohomology algebra $(H^{*, *}(M), \wedge)$ has some very
nice properties.
Let $\int_M: H^*(M, \bC) \rightarrow \bC$ be defined by the integrals of 
differential forms over $M$.
By Poincar\'{e} duality,
$\eta: H^{*, *}(M) \otimes H^{*, *}(M) \to \bC$ defined by
$$\eta(\alpha, \beta) = \int_M \alpha \wedge \beta$$
is a graded symmetric nondegenerate bilnear form on $H^{*, *}(M)$;
furthermore,
we have
$$\eta(\alpha \wedge \beta, \gamma) = \eta(\alpha, \beta \wedge \gamma)$$
for $\alpha, \beta, \gamma \in H^{*, *}(M)$.
In other words,
$(H^{*, *}(M), \wedge, \eta)$ is a bi-graded Frobenius algebra. 
By universal coefficient theorem, 
$H^*(M, \bQ) \otimes \bC \cong  H^*(M, \bC)$.
This shows that 
$(H^{*, *}(M), \wedge)$ is the complexification 
of a graded commutative algebra over $\bQ$.
The latter has a minimal model over $\bQ$,
which then determines a rational homotopy type.
Presumably,
it should be the same as 
the rational homotopy type of the K\"{a}hler manifold $M$.

Denote by $L$ the multiplication on forms by $\omega$ and $\Lambda$ its adjoint,
and $h = [\Lambda, L]$.
The commutation relations of these three operators are
that of $\fsl(2, \bR)$.
Notice that $L$ has bidegree $(1, 1)$,
$\Lambda$ has bidegree $(-1, -1)$
and $h$ has bidegree $(0, 0)$.
These operations induce a representation of $\fsl(2, \bR)$ 
on the space of harmonic forms 
which is used in the proof of the Hard Lefschetz Theorem attributed to Chern
(see e.g. Griffiths-Harris \cite{Gri-Har}).

\begin{definition}
A representation of $\fsl(2, \bR)$ on a bigraded vector space is
said to be {\em of Lefschetz type}
if $X$ acts by an operator of degree $(-1, -1)$,
$Y$ acts by an operator of bidegree $(1, 1)$
and $H$ acts by an operator of bidegree $(0, 0)$,
where $X$, $Y$, $H$ are generators of $\fsl(2, \bR)$
such that
\begin{align*}
[X, Y] &= H, & [H, X] &= 2X, & [H, Y] &= -2 Y.
\end{align*}
\end{definition}

To summarize, we have the following

\begin{theorem} \label{mainthm:Kahler}
For a closed K\"{a}hler manifold $M$,
$H^*(M, \bC)$ has a structure of finite dimensional
bi-graded Frobenius algebra which admits a representation 
by $\fsl(2, \bR)$ of Lefschetz type.   
Furthermore, $H^{*,*}(M)$ is the complexification 
of a rational graded commutative algebra..
\end{theorem}

\subsection{DGA $\cB$ of a Calabi-Yau manifold}
There is another natural DGA on a complex manifold:
$$\cB = \Omega^{-*, *}(M) = \oplus_{p, q} \Omega^{-p, q}(M) 
= \oplus_{p, q} \Gamma(M, \Lambda^pTM 
\otimes \Lambda^q\overline{T}^*M),$$
with wedge product $\wedge$, 
and differential 
$\bar{\partial} = \bar{\partial}_{\Lambda^*TM}$,
the $\bar{\partial}$ operator for the holomorphic vector bundle
$\Lambda^*TM$.
Elements of $\Omega^{-p, q}(M)$ are given the bidegree $(p, q)$.
The importance of this DGA is that it contains
the deformation complex of the complex structures:
$$\Omega^{0, 0}(TM) \stackrel{\bar{\partial}}{\rightarrow}
\Omega^{0, 1}(TM) \stackrel{\bar{\partial}}{\rightarrow}
\Omega^{0, 2}(TM),$$
hence we call $(\cB, \wedge, \bar{\partial})$ the {\em extended
deformation algebra}.
On a Calabi-Yau $n$-manifold $M$
with nontrivial holomorphic volume form $\Omega$, 
one can define an isomorphism 
$$\flat: \Omega^{-p, q}(M) \rightarrow \Omega^{n-p, q}(M),$$
and its inverse
$$\sharp: \Omega^{n-p, q}(M) \rightarrow \Omega^{-p, q}(M).$$
Tian \cite{Tia} and Todorov \cite{Tod}
defined an operator $\Delta$ of degree $-1$ on $\cB$ by
$$\Delta \gamma = (\partial \gamma^{\flat})^{\sharp}.$$
A formula in Tian \cite{Tia} shows that
$\Delta$ is related to the Schouten-Nijenhuis bracket $[\cdot, \cdot]$
on $\cB$ and reveals the DGBV algebra structure on
$(\cB, \wedge, \bar{\partial}, 
\Delta, [\cdot, \cdot])$.
This has very far-reaching consequences 
besides the original results of Tian and Todorov
on the deformations of Calabi-Yau manifold,
see Bershadsky-Cecotti-Ooguri-Vafa \cite{Ber-Cec-Oog-Vaf} and 
Barannikov-Kontsevich \cite{Bar-Kon}.
Following Tian and Todorov,
for any operator $\cO$ on $\cA$,
define an operator $\widetilde{\cO}$ on $\cB$ 
by setting $\widetilde{\cO}(\alpha) : = (\cO \alpha^{\flat})^{\sharp}$.
Then in this notation, $\Delta = \tilde{\partial}$,
$\sharp\square\flat = \widetilde{\square}$.
It is straightforward to see that
$\tilde{\bar{\partial}} = \bar{\partial}_{\Lambda^*TM}$.
It is clear that for two operators $\cO_1$ and $\cO_2$ on $\cA$,
we have $\widetilde{\cO_1\cO_2} = \widetilde{\cO_1}\widetilde{\cO_2}$,
hence $\widetilde{[\cO_1, \cO_2]} = [\widetilde{\cO_1}, \widetilde{\cO_2}]$.
Since $\flat$ and $\sharp$ are isometries, it is easy to see that 
$\widetilde{\cO^*} = \widetilde{\cO}^*$ for any operator $\cO$ on $\cA$.
From (\ref{eqn:comm1}) to (\ref{eqn:comm4}),
we deduce the following identities:
\begin{align*}
& [\tilde{\partial}, \tilde{\partial}] = 0, 
	& [\tilde{\bar{\partial}}, \tilde{\bar{\partial}}] &= 0, 
	& [\tilde{\partial}, \tilde{\bar{\partial}}] &= 0, \\
& [\tilde{\partial}^*, \tilde{\partial}^*] = 0, 
	& [\tilde{\bar{\partial}}^*, \tilde{\bar{\partial}}^*] &= 0, 
	& [\tilde{\partial}^*, \tilde{\bar{\partial}}^*] &= 0, \\
& [\tilde{\partial}, \tilde{\bar{\partial}}^*] = 0, 
	& [\tilde{\partial}^*, \tilde{\bar{\partial}}] &= 0,
	& \\
& [\tilde{\partial}, \tilde{\partial}^*] = \frac{1}{2}\widetilde{\square}, 
& [\tilde{\bar{\partial}}, \tilde{\bar{\partial}}^*] 
	& = \frac{1}{2}\widetilde{\square}, 
	& 
\end{align*}
As usual, we rewrite these operators 
in terms of the Levi-Civita connection $\nabla$.
For any $x \in M$, let $\{e^1_x, \cdots, e^n_x\}$ be a basic 
of the complex cotangent space $T^*M$,
such that 
$\omega_x = e^1_x \wedge \bar{e}^1 + \cdots e^n_x \wedge \bar{e}^n_x$.
Extend this basis to a local frame $\{e^1, \cdots, e^n\}$ 
in a neighborhood $U$ of $x$ by parallel transportations
along the geodesics starting from $x$.
Then we have 
$\omega = e^1 \wedge \bar{e}^1 + \cdots e^n \wedge \bar{e}^n$ in $U$.
Let $\{e_i, \cdots, e_n\}$ be the dual frame for $TM$, 
then we have   
$$\nabla_{e_i} e^j = \nabla_{\bar{e}_i} e^j 
= \nabla_{e_i} \bar{e}^j =\nabla_{\bar{e}_i} \bar{e}^j = 0$$
at $x$. 
On $U$, we have
\begin{align*}
& \partial = e^i \wedge \nabla_{e_i}, 
& \bar{\partial} & = \bar{e}^i \wedge \nabla_{\bar{e}_i}, \\
& \partial^* = -e_i \vdash \nabla_{\bar{e}_i},
& \bar{\partial}^* & = - \bar{e_i} \vdash\nabla_{e_i},
\end{align*}
where $\vdash$ stands for contractions.
Up to a constant, we have $\Omega = e^1 \wedge \cdots \wedge e^n$.
Since $\Omega$ is parallel, 
$\nabla_{e_i}$ and $\nabla_{\bar{e}_i}$ commute with $\sharp$ and $\flat$.
Then one sees that
\begin{align*}
& \tilde{\partial} = e^i \vdash \nabla_{e_i}, 
& \tilde{\bar{\partial}} & = (-1)^n \bar{e}^i \wedge \nabla_{\bar{e}_i}, \\
& \tilde{\partial}^* = -e_i \wedge \nabla_{\bar{e}_i},
& \tilde{\bar{\partial}}^* & = -(-1)^n \bar{e}_i \vdash\nabla_{e_i}.
\end{align*}
Now it is clear that $\tilde{\partial}^*$ and $\tilde{\bar{\partial}}$ 
are derivations of degrees $(1, 0)$ and $(0, 1)$ respectively.
They satisfy the conditions in Theorem \ref{thm:formal},
so we get
 
\begin{theorem} \label{thm:CY}
For a Calabi-Yau manifold $M$, 
the DGA $(\cB, \wedge, \bar{\partial})$ is formal.
\end{theorem}

The algebra $(H^{-*, *}(M), \wedge)$ shares
some properties with $(H^{*, *}(M), \wedge)$.
Barannikov and Kontsevich used an integral 
$\int: \Omega^{-*, *}(M) \rightarrow \bC$ by
$$\int \gamma : = \int_M \gamma^{\flat} \wedge \Omega.$$
It satisfies 
$$\int \bar{\partial} \alpha \wedge \beta =
(-1)^{|\alpha|+1} \int \alpha \wedge \bar{\partial}\beta.$$
(See Claim 4.1 in \cite{Bar-Kon}. Note that we use a 
different grading here.) 
It induces a structure of a bigraded Frobenius algebra on $H^{-*, *}(M)$.
For a general Calabi-Yau manifold,
mirror symmetry suggests the existence of an element 
$\tilde{\omega} \in H^{-1, 1}(M)$ which plays the 
role of the K\"{a}hler form $\omega$, 
in the sense that the three operators given by
the multiplication by $\tilde{\omega}$, its adjoint,
and their commutator
yield a representation of $\fsl(2, \bR)$ on
$H^{-*, *}(M)$ of Lefschetz type.
However, one should not expect that $\tilde{\omega}$ 
is determined by $\omega$ and $\Omega$,
except for hyperk\"{a}hler manifolds,
so extra work is required to establish its existence.
For Calabi-Yau hypersurfaces in weighted projective spaces,
 H\"{u}bsch and Yau \cite{Hub-Yau}
have found such an element and defined an $SL(2, \bC)$ 
representation on $\oplus_q H^{-q, q}(M)$.
F general Calabi-Yau manifolds,
some preliminary results have been obtained by 
Cao and the author \cite{Cao-Zho5}. 
On the other hand a non-Lefschetz type representation 
can be easily found as follows:
consider operators $\tilde{L}$, $\tilde{\Lambda}$ and $\tilde{h}$.
They form an $\fsl(2, \bR)$ algebra. 
However $\tilde{L}$ and $\tilde{\Lambda}$ have bidegrees
$(-1, 1)$ and $(1, -1)$ respectively.
The usual commutation relations among 
$L, \Lambda, h, \partial, \bar{\partial}, \partial^*, 
\bar{\partial}^*, \square$ translate to
the commutation relations among
$\tilde{L}, \tilde{\Lambda}, \tilde{h}$ etc.
One then gets a representation of $\fsl(2, \bR)$ on $H^{-*, *}(M)$.

Use the isomorphisms $\flat$ and $\sharp$,
the study of $H^{-*, *}(M)$ can be transformed to the 
study of $H^{*, *}(M)$.
The wedge produce $\wedge$ on $H^{-*, *}(M)$  is transformed 
to product $\tilde{\wedge}$ on $H^{*, *}(M)$ defined by
$$\alpha \tilde{\wedge} \beta = 
(\alpha ^{\sharp} \wedge \beta^{\sharp})^{\flat}.$$
The integral $\int$ on $H^{-*, *}(M)$ is
transformed to the integral $\tilde{\int}_M$on $H^{*, *}(M)$ given by
$$\tilde{\int}_M \alpha = \int_M \alpha \wedge \Omega.$$
Therefore, we get a Frobenius algebra $(H^{*, *}(M), \tilde{\wedge})$
which is isomorphic to the Frobenius algebra $(H^{-*, *}(M), \wedge)$.
Now let $\{ \gamma_a \}$ be a basis of $H^{*, *}(M)$,
such that each $\gamma_a$ represents a class 
which lies in the image of $H^*(M, \bQ) \rightarrow H^*(M, \bC)$,
and $\gamma_1 = 1$.
We say the holomorphic volume form $\Omega$
satisfies the rationality condition
if 
$$\Phi_{abc} = \int (\gamma_a^{\sharp}\wedge \gamma_b^{\sharp}
\wedge \gamma_c^{\sharp})^{\flat} \wedge \Omega$$
are all rational.
When this condition holds,
the graded commutative algebra $(H^{-*, *}(M), \wedge)$ is
the complexification of a rational graded commutative algebra.
When $M$ is simply connected,
the degree $1$ part of $H^{-*, *}(M)$ is trivial:
$H^{0, 1}(M)$ is clearly trivial,
and $H^{-1, 0} \cong H^{n-1, 0}(M) = 0$.
By rational homotopy theory,
we have the following 

\begin{theorem} \label{mainthm:CY}
If the holomorphic volume form $\Omega$ on a simply connected
Calabi-Yau manifold $M$
satisfies the rationality condition,
then the graded commutative algebra $(H^{-*, *}(M), \wedge)$
 determines a unique rational homotopy type.
\end{theorem}

Under the condition of the above theorem , 
if $M$ has a mirror manifold,
i.e. a K\"{a}hler manifold  $\widetilde{M}$ such that  
$(H^{*, *}(\widetilde{M}), \wedge) \cong (H^{-*,*}(M), \wedge)$
as bigraded algebras,
byTheorem \ref{thm:formal2},
$\widetilde{M}$ should have the rational homotopy determined by
$(H^{-*,*}(M), \wedge)$.
On the other hand,
it is possible that the rationality condition 
is necessary for the existence of a mirror manifold.

If $M$ be a hyperk\"{a}hler manifold
then $M$ has a holomorphic symplectic form which induces an
isomorphism of DGA's
$$(\Omega^{-*, *}(M), \wedge, \bar{\partial}) \rightarrow 
(\Omega^{*, *}(M), \wedge, \bar{\partial}).$$
Therefore, a closed hyperk\"{a}hler manifold 
should be the mirror manifold of itself.

\section{Further discussions}

We have considered a single Calabi-Yau manifold so far. 
It is important in the study of mirror symmetry to consider the
moduli space of complex structures and K\"{a}hler strctures,
then one will encounter differentiably different Calabi-Yau manifolds
which lies in different but connected regions of such moduli spaces. 
It is very likely that these different Calabi-Yau manifolds 
determines the same complex or rational homotopy types by their
Dolbeault or deformation algebras. 
For moduli space of complex structure structures,
this is easy to understand.
A small deformation of a complex manifold will not change its 
topological type.
However,
a large deformation may lead to a singularity and then a birationally related
Calabi-Yau manifold.
Now birationally related Calabi-Yau manifolds can be obtained from
each other by flops.
It might be possible to directly check that flops do 
not change the rational homotopy type.
For example,
Tian \cite{Tia} studied the smoothing of Calabi-Yau 3-folds with ordinary
double points.
It is possible that such singular 3-folds lies on some walls that
divide Calabi-Yau 3-folds related by flops.
Earlier, 
Tian and Yau \cite{Tia-Yau} gave examples of Calabi-Yau manifolds which
are related to each other by flops. 
They are different as differential manifolds, 
but they have the same cohomology ring,
hence have the same complex homotopy types. 
Similarly,
one can glue the K\"{a}hler moduli spaces of
 birationally related Calabi-Yau manifolds or orbifolds
along the degenerations of the K\"{a}hler structure on
their complex subvarieties 
to form the enlarged K\"{a}hler moduli space.
See e.g. Aspinwall-Greene-Morrison \cite{Asp-Gre-Mor}.
Again the relevant Calabi-Yau 3-folds 
should have the same rational homotopy type.
In the above discussion,
we expect the complex homotopy types determined by
the deformation algebras are also the same.

\end{document}